\documentclass[a4paper,pdftex,reqno,10pt]{article}

\usepackage{geometry}
\usepackage{graphicx}
\usepackage{multirow}
\usepackage{times}
\usepackage{amssymb,amsmath}

\newcommand{\bbox}{\vrule height7pt width4pt depth1pt}
\def\@begintheorem#1#2{\list{}{\thm@body}%
  \item[]{\bf #1~#2.}\quad\it\ignorespaces}
\def\@opargbegintheorem#1#2#3{\list{}{\thm@body}%
  \item[]{\bf #1~#2~\ifrembrks #3\global\rembrksfalse\else (#3)\fi.}%
  \quad\it\ignorespaces}
\def\@endtheorem{\endlist}
\newtheorem{theorem}{Theorem}[section]

\newtheorem{lemma}[theorem]{Lemma}

\begin{document}

\title{Regular matchstick graphs}
\author{
Sascha Kurz\thanks{sascha.kurz@uni-bayreuth.de}\\
Fakult\"at f\"ur Mathematik, Physik und Informatik, Universit\"at Bayreuth, Germany
\and
Rom Pinchasi\thanks{room@math.technion.ac.il}\\
Mathematics Dept., Technion---Israel Institute of Technology, Haifa 32000, Israel
}

\maketitle

\noindent
  \textbf{Abstract.}
  A graph $G=(V,E)$ is called a unit-distance graph in the plane if there is an injective embedding of $V$ in
  the plane such that every pair of adjacent vertices are at unit distance apart. If additionally the corresponding
  edges are non-crossing and all vertices have the same degree $r$ we talk of a regular matchstick graph. Due to
  Euler's polyhedron formula we have $r\le 5$. The smallest known $4$-regular matchstick graph is the so called
  Harborth graph consisting of $52$ vertices. In this article we prove that no finite $5$-regular matchstick graph
  exists and provide a lower bound for the number of vertices of $4$-regular matchstick graphs.

\noindent
{
  \center\small{Keywords: unit-distance graphs} 
}
\noindent
\rule{\textwidth}{0.3 mm}

\section{Introduction}
\noindent
One of the possibly best known problems in combinatorial geometry asks how often the same distance can occur among
$n$ points in the plane. Via scaling we can assume that the most frequent distance has length $1$. Given any set $P$ of points in the plane, we can define the so called unit-distance graph in the plane, connecting two elements of $P$ by an edge if their distance is one. The known bounds for the maximum number $u(n)$ of edges of a unit-distance graph in the plane, see e.~g.{} \cite{1086.52001}, are given by
\[
  \Omega\!\left(ne^{\frac{c\log n}{\log\log n}}\right)\le u(n)\le O\!\left(n^{\frac{4}{3}}\right).
\]
For $n\le 14$ the exact numbers of $u(n)$ were determined in \cite{schade}, see also \cite{1086.52001}.

If we additionally require that the edges are non-crossing, then we obtain another class of geometrical and combinatorial objects:
A \textbf{matchstick graph} is graph drawn with straight edges in the plane  such that the edges have unit length, and non-adjacent edges do not intersect, see Figure~\ref{fig_configuration_1} for an example.
%

For matchstick graphs the known bounds for the maximum number $\tilde{u}(n)$ of edges, see e.~g.{} \cite{1086.52001}, are given by
\[
  \left\lfloor 3n-\sqrt{12n-3}\right\rfloor\le \tilde{u}(n)\le 3n-O\!\left(\sqrt{n}\right),
\]
where the lower bound is conjectured to be exact.

We call a matchstick graph $r$-regular if every vertex has degree $r$. In \cite{matchsticks_in_the_plane} the authors consider $r$-regular matchstick graphs with the minimum number $m(r)$ of vertices. Obviously we have $m(0)=1$, $m(1)=2$, and $m(2)=3$, corresponding to a single vertex, a single edge, and a triangle, respectively. The determination of $m(3)=8$ is an amusement left to the reader. 
For degree $r=4$ the exact determination of $m(4)$ is unsettled so far. The smallest known example is the so called Harborth graph, see e.~g.{} \cite{gerbracht} for a drawing, yielding $m(4)\le 52$. Here we prove $m(4)\ge 20$.

Due to the Eulerian polyhedron formula every finite planar graph contains a vertex of degree at most five so that we have $m(r)=\infty$ for $r\ge 6$. Examples are given by the regular triangular lattice or an infinite $6$-regular tree.

For degree $5$ it is announced at several places that no finite $5$-regular matchstick graph does exist. 
%
%
In \cite{1063.05036} the authors list five publications where they have mentioned an unpublished proof\footnote{In the meantime the temporally lost manuscript was recovered and after some minor corrections the arguments turned out to be valid. A retyped and slightly edited electronic version can be found at http://www.wm.uni-bayreuth.de/fileadmin/Sascha/aartb.pdf.} for the non-existence of a $5$-regular matchstick graph and state that up to their knowledge the problem is open. Currently there is no published proof available.\footnote{There is a further unpublished manuscript containing a rather technical and long proof \cite{regular_matchstick2}.}
%

The aim of this article is indeed to bridge this gap.



\section{$5$-regular matchstick graphs}
\begin{theorem}
  \label{thm_main}
  No finite $5$-regular matchstick graph does exist.
\end{theorem}

  \textsc{Proof.}
  Suppose to the contrary that there is such a graph and let $M$ be the planar
  map which is drawing of this graph in the plane such that every edge is a unit 
  length straight line segment and no two edges cross.

  Without loss of generality we assume that this graph is connected and denote
  by $V$ the number of its vertices, by $E$ the number of its edges, and by $F$
  the number of faces in the planar map $M$. By Euler's formula we have
  $V-E+F=2$. For every $k \geq 3$ we denote by $f_{k}$ the number of faces
  in $M$ with precisely $k$ edges.

  We observe that $2E=\sum kf_{k}=5V$ and $F=\sum f_{k}$. Therefore,
  \begin{equation}\label{eq:1}
    -6=
    -3V+E+2E-3F=-3V+\frac{5}{2}V+\sum kf_{k}-3\sum f_{k}=
    -\frac{1}{2}V+\sum(k-3)f_{k}.
  \end{equation}
  We begin by giving a charge of $-\frac{1}{2}$ to each vertex and by giving
  a charge of $k-3$ to each face in $M$ with precisely $k$ edges.
  By~(\ref{eq:1}) the total charge of all the vertices and faces is negative.
  We will reach a contradiction by redistributing the charge in such a way
  that eventually each vertex and each face will have a non-negative charge.

  We redistribute the charge in the following very simple way.
  Consider a face $T$ of $M$ and a vertex $x$ of $T$. Let $\alpha$ denote 
  measure of the internal angle of $T$ at $x$. Only if $\alpha > \frac{\pi}{3}$
  we take a charge of 
  $\min\!\left(\frac{1}{2}, \frac{3}{2\pi}\alpha-\frac{1}{2}\right)$ from $T$ and move it to 
  $x$.

  We now show that after the redistribution of charges every vertex and every 
  face have a non-negative charge. Consider a vertex $x$. Let $\ell$ denote the
  number of internal angles at $x$ that are greater than $\frac{\pi}{3}$.
  As the degree of $x$ equals to $5$ we must have $\ell >0$.
  If due to one of these $\ell$ angles we transfered a charge of $\frac{1}{2}$
  to $x$, then the charge at $x$ is non-negative. Otherwise, note the the sum of 
  these $\ell$ angles is at least $2\pi-(5-\ell)\frac{\pi}{3}=
  \frac{\pi}{3}(\ell+1)$. Hence the total charge transfered to $x$ due to these 
  angles is at least $\frac{3}{2\pi}\frac{\pi}{3}(\ell+1)-\frac{\ell}{2}=
  \frac{1}{2}$. Here again we conclude that the charge at $x$ is non-negative.

  Consider now a face $T$ in $M$ with $k \geq 3$ edges. Assume first that $T$ is 
  a bounded face. The initial charge of $T$ is $k-3 \geq 0$. Therefore, 
  if the charge at $T$ becomes negative this implies that one of the internal
  angles of $T$ is greater than $\frac{\pi}{3}$. In particular $T$ cannot be a 
  triangle and thus $k \geq 4$. If $k=4$, then $T$ is a rhombus. If only two 
  internal angles of $T$ are greater than $\frac{\pi}{3}$, then at most
  a total charge of $1$ was deduced from the initial charge of $T$, leaving
  its charge non-negative. If all internal angles of $T$ are greater than 
  $\frac{\pi}{3}$, then the total charge deduced from $T$ is at most
  $\frac{3}{2\pi}\cdot 2\pi-\frac{4}{2}=1$, leaving the charge at $T$ non-negative. 

  If $k=5$ and the charge of $T$ is negative after we
  redistribute the charges, then each of the internal angles of $T$ must be 
  greater than $\frac{\pi}{3}$ and as the sum of the internal angles of $T$ is
  equal to $3\pi$, the charge deduced from $T$ amounts to at most
  $\frac{3}{2\pi}\cdot 3\pi-\frac{5}{2}=2$, leaving the charge at $T$ 
  non-negative. 

  Finally if $k \geq 6$, then the charge deuced from $T$ is at most
  $\frac{k}{2}$ leaving a charge of at least $k-3-\frac{k}{2} \geq 0$.

  It is left to consider the unbounded face $S$ of $M$. If the number of edges
  of $S$ is at least $6$, we are done as in the case of a bounded face. 
  The cases where the unbounded 
  face consist of at most $5$ edges can be easily excluded. Another way to settle
  this issue is to observe that if $S$ consists of at most $5$ edges, then 
  the total charge deduced from $S$ is at most $\frac{5}{2}$ leaving the charge 
  of $S$ at least $-\frac{5}{2}$ (and in fact at least $-\frac{3}{2}$). 
  We still obtain a contradiction as the sum of all charges should be equal to
  $-6$ while only the unbounded face may remain with a negative charge that
  is not smaller than $-\frac{5}{2}$. \hfill{\bbox}

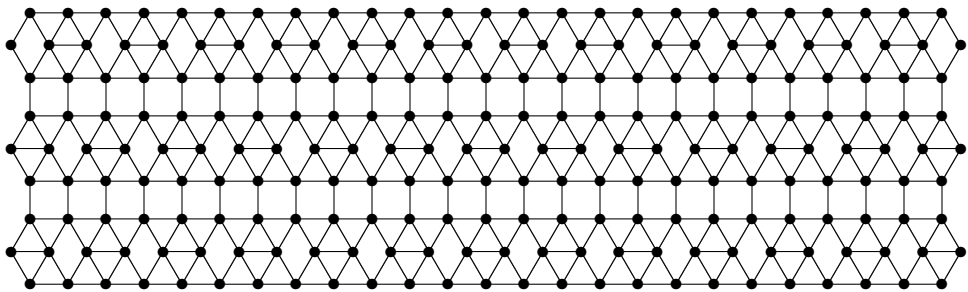
\begin{figure}[ht]
  \begin{center}
    \setlength{\unitlength}{0.5cm}
    \begin{picture}(26,7.2)
      \multiput(0,0)(1,0){25}{\qbezier(0,0.866)(0.25,0.433)(0.5,0)}
      \multiput(0,0)(1,0){25}{\qbezier(0,0.866)(0.25,1.299)(0.5,1.732)}
      \multiput(0,0)(1,0){25}{\qbezier(1,0.866)(0.75,0.433)(0.5,0)}
      \multiput(0,0)(1,0){25}{\qbezier(1,0.866)(0.75,1.299)(0.5,1.732)}
      \multiput(0,0.866)(1,0){26}{\circle*{0.3}}
      \multiput(0.5,0)(1,0){25}{\circle*{0.3}}
      \multiput(0.5,1.732)(1,0){25}{\circle*{0.3}}
      \put(0.5,0){\line(1,0){24}}
      \put(0.5,1.732){\line(1,0){24}}
      \multiput(0,0.866)(2,0){13}{\line(1,0){1}}
      \multiput(0.5,1.732)(1,0){25}{\line(0,1){1}}
      \multiput(0,2.732)(1,0){25}{\qbezier(0,0.866)(0.25,0.433)(0.5,0)}
      \multiput(0,2.732)(1,0){25}{\qbezier(0,0.866)(0.25,1.299)(0.5,1.732)}
      \multiput(0,2.732)(1,0){25}{\qbezier(1,0.866)(0.75,0.433)(0.5,0)}
      \multiput(0,2.732)(1,0){25}{\qbezier(1,0.866)(0.75,1.299)(0.5,1.732)}
      \multiput(0,3.598)(1,0){26}{\circle*{0.3}}
      \multiput(0.5,2.732)(1,0){25}{\circle*{0.3}}
      \multiput(0.5,4.464)(1,0){25}{\circle*{0.3}}
      \put(0.5,2.732){\line(1,0){24}}
      \put(0.5,4.464){\line(1,0){24}}
      \multiput(0,3.598)(2,0){13}{\line(1,0){1}}
      \multiput(0.5,4.464)(1,0){25}{\line(0,1){1}}
      \multiput(0,5.464)(1,0){25}{\qbezier(0,0.866)(0.25,0.433)(0.5,0)}
      \multiput(0,5.464)(1,0){25}{\qbezier(0,0.866)(0.25,1.299)(0.5,1.732)}
      \multiput(0,5.464)(1,0){25}{\qbezier(1,0.866)(0.75,0.433)(0.5,0)}
      \multiput(0,5.464)(1,0){25}{\qbezier(1,0.866)(0.75,1.299)(0.5,1.732)}
      \multiput(0,6.340)(1,0){26}{\circle*{0.3}}
      \multiput(0.5,5.464)(1,0){25}{\circle*{0.3}}
      \multiput(0.5,7.196)(1,0){25}{\circle*{0.3}}
      \put(0.5,5.464){\line(1,0){24}}
      \put(0.5,7.196){\line(1,0){24}}
      \multiput(1,6.340)(2,0){12}{\line(1,0){1}}
    \end{picture}
  \end{center}
  \caption{Infinite $5$-regular match stick graph.}
  \label{fig_infinite}
\end{figure}

In Figure \ref{fig_infinite} we have drawn a fraction of an infinite $5$-regular matchstick graph. Since we can shift single rows of such a construction there exists an uncountable number of these graphs (even in the combinatorial sense). As mentioned by Bojan Mohar there are several further examples. 

Starting with the infinite $6$-regular triangulation we can delete several vertices at different rows and obtain an example by suitably removing horizontal edges to enforce degree $5$ for all vertices. Taking an arbitrary planar unit-distance graph with maximum degree $5$, whose vertices not being part of the unbounded face have degree exactly $5$, and continuing the  boundary vertices with trees gives another set of examples.

\section{$4$-regular matchstick graphs}

\noindent
In this section we show that every $4$-regular matchstick graph must 
consist of at least $20$ vertices. Although this is a little far from the
best known construction with $52$ vertices (\cite{1063.05036,gerbracht,matchsticks_in_the_plane}),
we bring some arguments that may lead the way to obtain improved bounds and
remark, that the whole proof is free from computer calculations and massive
case differentiations. For a determination of the lower bound $m(4)\ge 34$,
based on massive computer calculations, we refer the interested reader to 
\cite{unmasking}.

\begin{lemma}
  \label{lemma_bound_4_regular}
  For every connected $4$-regular matchstick graph in the plane we have
  $
    n\ge \frac{3k-l+12}{2}
  $,
  where $n$ denotes the number of vertices, $k$ the number of edges of the unbounded face, and
  $l\le k-3$ the number of internal angles of the unbounded face that are equal to $\pi$
\end{lemma}
  \textsc{Proof.}
  Assume $G$ is a $4$-regular matchstick graph.
  Let us denote by $V$ the number of its vertices, by $E$ the number of its edges, and by $F$ the number
  of its faces in the planar map $M$ realizing $G$ as a matchstick graph in 
  the plane. By Euler's formula we have $V-E+F=2$. For every $r\ge 3$ we denote
  by $f_r$ the number of faces in $M$ with precisely $r$ edges.

  We observe that $2E=\sum rf_r=4V$ and $F=\sum f_r$. Therefore,
  \begin{eqnarray}
    -6 & = & -3V+3E-3F=-3V+E+2E-3F=\nonumber\\
    & = & -3V+\frac{4}{2}V+\sum rf_r-3\sum f_r=-V+\sum(r-3)f_r.\label{eq_2}
  \end{eqnarray}
  We begin by giving a charge of $-1$ to each vertex and by giving a charge of $r-3$ to each face in $M$
  with precisely $r$ edges. By~(\ref{eq_2}) the total charge of all the vertices and faces is $-6$.
  
  We now redistribute the charge in the following simple way. Consider a face $T$ of $M$ and a vertex $x$ of $T$.
  Let $\alpha$ denote the measure of the internal angle of $T$ at $x$. If $T$ is an equilateral triangle, we do
  nothing. If $\alpha=\pi$ we move a charge of $\frac{7}{8}$ from $T$ to $x$. Otherwise, we take a charge of
  $\frac{9}{8\pi}\alpha-\frac{5}{16}$ from $T$ and move it to $x$.
  
  We now show that after the redistribution of charges every vertex is at least $-\frac{1}{8}$ and every bounded face has a
  nonnegative charge.

  For the first part consider a vertex $x$. Let $p$ denote the number of internal angles at $x$ which are equal to
  $\frac{\pi}{3}$. Clearly, $0\le p\le 3$. If $l=3$, then one of the internal angles at $x$ is equal to $\pi$. It follows
  that the charge at $x$ is at least $-1+\frac{7}{8}$. This is because $x$ never gets a negative contribution from
  internal angles that are equal to $\frac{\pi}{3}$.
  
  If $l\le 2$, then the charge at $x$ is at least
  $$
    -1+\frac{9}{8\pi}\left(2\pi-\frac{l\pi}{3}\right)-\frac{5(4-l)}{16}=-\frac{l}{16}\ge-\frac{1}{8}.
  $$
  Consider now a bounded face $T$ in $M$ with $r\ge 3$ edges. We will show that the charge of $T$ after we have
  redistributed the charges is nonnegative. The initial charge of $T$ is $r-3\ge 0$. If $T$ is a (necessarily
  equilateral) triangle, then the charge of $T$ does not change and it remains $0$.
  If $T$ is a quadrilateral and hence a rhombus, then clearly no internal angle of $T$ can be equal to $\pi$ and the
  total charge removed from $T$ is $\frac{9}{8\pi}\cdot 2\pi-\frac{5}{16}\cdot 4=1$. This leaves $T$ with a total charge
  of $0$.

  Suppose now that $T$ has $r\ge 5$ edges. The sum of all internal edges of $T$ is equal to $\pi(r-2)$. Let $p$ denote
  the number of internal angles of $T$ that are equal to $\pi$. Clearly, $p\le r-3$. The total charge removed from $T$
  is equal to
  $$
    \frac{7}{8}p+\frac{9}{8\pi}\cdot\pi(r-2-p)-\frac{5(r-p)}{16}=\frac{p}{16}+\frac{13}{16}r-\frac{9}{4}.
  $$
  Because $p\le r-3$ this charge is at most
  $$
    \frac{r-3}{16}+\frac{13}{16}r-\frac{9}{4}=\frac{7}{8}r-\frac{39}{16}.
  $$
  This is smaller than the initial charge of $T$ which is $r-3$ for all $r\ge 5$.

  We have thus shown that every bounded face of $M$ remains with a nonnegative charge. But what about the unbounded
  face. Suppose $T$ is the unbounded face and it has $k$ edges. The sum of all internal edges of $T$ is equal to
  $\pi(k+2)$. Let $l$ denote the number of internal angles of $T$ that are equal to $\pi$. Here too we have $l\le k-3$.
  
  The total charge removed from $T$ is equal to
  $$
    \frac{7}{8}l+\frac{9}{8\pi}\cdot \pi(k+2-l)-\frac{5(k-l)}{16}=\frac{l}{16}+\frac{13}{16}k+\frac{9}{4}.
  $$
  This leaves $T$ with a total charge of $\frac{3}{16}k-\frac{l}{16}-\frac{21}{4}$.

  The total charge of all vertices and faces should be equal to $-6$. The total charge of all vertices is at least
  $-\frac{n}{8}$ and every bounded face has a nonnegative charge. Therefore, we have
  $$
    \frac{3}{16}k-\frac{l}{16}-\frac{21}{4} -\frac{n}{8} \le -6,
  $$
  which is equivalent to
  $$
    n\ge \frac{3k-l+12}{2}.
  $$\hfill{\bbox}

From $2E=\sum rf_r=4V$, $F=\sum f_r$, and Euler's formula $V-E+F=2$ we conclude
\begin{equation}
  \sum_{r=3}^\infty (4-r)f_r=f_3-f_5-2f_6-3f_7-4f_8-\dots =8\label{eq_A_i_sum}, 
\end{equation}
using the notation from the proof of Lemma \ref{lemma_bound_4_regular}. Thus we have $f_3\ge k+4$,
where $k$ denote the number of edges of the unbounded face. The maximum area of an equilateral $k$-gon,
whose edges have length $1$, is given by
\begin{equation}
  A_{\text{max}}(k)=\frac{k}{4}\cdot\cot\left(\frac{\pi}{k}\right).
\end{equation}
Thus $k$ must be big enough in order for the complement of the unbounded face 
to contain at least $k+4$ triangles. In Table~\ref{table_max_area} we have
listed the maximum area of an equilateral $k$-gon measured in units of equilateral triangles.

\begin{table}[ht]
  \begin{tabular}{rrrrrrrrrrr}
    \hline
    $k$ & 3 & 4 & 5 & 6 & 7 & 8 & 9 & 10 & 11 & 12 \\
    $\frac{A_{\text{max}}(k)}{A_{\text{max}}(3)}$ &
    1.000 & 2.310 & 3.974 & 6.000 & 8.393 & 11.151 & 14.277 & 17.770 & 21.630 & 25.857 \\
    \hline
  \end{tabular}
  \caption{Maximum number of equilateral triangles inside an equilateral $k$-gon.}
  \label{table_max_area}
\end{table}

Solving $\frac{k}{4}\cdot\cot\left(\frac{\pi}{k}\right)\ge (k+4)\cdot\frac{\sqrt{3}}{4}$ or using Table~\ref{table_max_area}
yields $k\ge 9$ which can be improved by some further arguments.

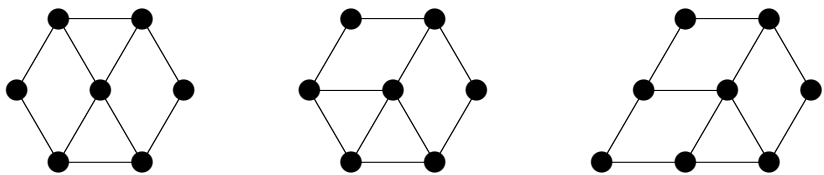
\begin{figure}[h]
  \begin{center}
    \setlength{\unitlength}{1.10cm}
    \begin{picture}(9,1.8)
      \put(0,0.866){\circle*{0.25}}
      \put(0.5,0){\circle*{0.25}}
      \put(0.5,1.732){\circle*{0.25}}
      \put(1,0.866){\circle*{0.25}}
      \put(1.5,0){\circle*{0.25}}
      \put(1.5,1.732){\circle*{0.25}}
      \put(2,0.866){\circle*{0.25}}
      \put(0.5,0){\line(1,0){1}}
      \put(0.5,1.732){\line(1,0){1}}
      \qbezier(0,0.866)(0.25,0.433)(0.5,0)
      \qbezier(0,0.866)(0.25,1.299)(0.5,1.732)
      \qbezier(1,0.866)(0.75,0.433)(0.5,0)
      \qbezier(1,0.866)(0.75,1.299)(0.5,1.732)
      \qbezier(1,0.866)(1.25,0.433)(1.5,0)
      \qbezier(1,0.866)(1.25,1.299)(1.5,1.732)
      \qbezier(2,0.866)(1.75,0.433)(1.5,0)
      \qbezier(2,0.866)(1.75,1.299)(1.5,1.732)
      \put(3.5,0.866){\circle*{0.25}}
      \put(4,0){\circle*{0.25}}
      \put(4,1.732){\circle*{0.25}}
      \put(4.5,.866){\circle*{0.25}}
      \put(5,0){\circle*{0.25}}
      \put(5,1.732){\circle*{0.25}}
      \put(5.5,0.866){\circle*{0.25}}
      \put(3.5,0.866){\line(1,0){1}}
      \put(4,0){\line(1,0){1}}
      \put(4,1.732){\line(1,0){1}}
      \qbezier(3.5,0.866)(3.75,0.433)(4,0)
      \qbezier(3.5,0.866)(3.75,1.299)(4,1.732)
      \qbezier(4.5,0.866)(4.25,0.433)(4,0)
      \qbezier(4.5,0.866)(4.75,0.433)(5,0)
      \qbezier(4.5,0.866)(4.75,1.299)(5,1.732)
      \qbezier(5.5,0.866)(5.25,0.433)(5,0)
      \qbezier(5.5,0.866)(5.25,1.299)(5,1.732)
      \put(7.5,0.866){\circle*{0.25}}
      \put(8,0){\circle*{0.25}}
      \put(7,0){\circle*{0.25}}
      \put(8,1.732){\circle*{0.25}}
      \put(8.5,.866){\circle*{0.25}}
      \put(9,0){\circle*{0.25}}
      \put(9,1.732){\circle*{0.25}}
      \put(9.5,0.866){\circle*{0.25}}
      \put(7.5,0.866){\line(1,0){1}}
      \put(7,0){\line(1,0){1}}
      \put(8,0){\line(1,0){1}}
      \put(8,1.732){\line(1,0){1}}
      \qbezier(7.5,0.866)(7.25,0.433)(7,0)
      \qbezier(7.5,0.866)(7.75,1.299)(8,1.732)
      \qbezier(8.5,0.866)(8.25,0.433)(8,0)
      \qbezier(8.5,0.866)(8.75,0.433)(9,0)
      \qbezier(8.5,0.866)(8.75,1.299)(9,1.732)
      \qbezier(9.5,0.866)(9.25,0.433)(9,0)
      \qbezier(9.5,0.866)(9.25,1.299)(9,1.732)
    \end{picture}\\[2mm]
    \caption{Configurations of triangles and quadrangles.}
    \label{fig_configuration_1}
  \end{center}
\end{figure}

\begin{lemma}
  \label{lemma_inner_3_4_configuration}
  If a $4$-regular matchstick graph contains an inner vertex but no bounded 
  faces with $s\ge 5$ edges, then it
  contains one of the configurations of Figure \ref{fig_configuration_1}.
\end{lemma}
  \textsc{Proof.}
  Due to the premises all inner faces are triangles and quadrangles. An inner vertex is part of four inner
  faces. Due to angle sums at most two of them can be triangles.\hfill{\bbox}

\begin{lemma}
  \label{lemma_inner_configuration_1}
  The total area of the quadrangles in each of the three configurations of Figure~\ref{fig_configuration_1},
  where the sides have length $1$, is greater than $\frac{\sqrt{3}}{2}$.
\end{lemma}
  \textsc{Proof.}
  If we denote the angles of the quadrangles at the central vertex by $\alpha_i$ the 
  area of the quadrangles is given by $\sum_i\left|\sin\alpha_i\right|$. A little calculus using 
  $0\le\alpha_i\le \pi$ yields the stated result.\hfill{\bbox}

\begin{lemma}
  \label{lemma_area_odd}
  For integers $s\ge 1$ the minimum area of an equilateral $(2s+1)$-gon with side lengths $1$ is at least $\frac{\sqrt{3}}{4}$.
\end{lemma}
  \textsc{Proof.}
  See \cite{0968.51014}.\hfill{\bbox}

Thus we can conclude $\frac{k}{4}\cdot\cot\left(\frac{\pi}{k}\right)\ge (k+6)\cdot\frac{\sqrt{3}}{4}$ yielding $k\ge 10$.

\begin{lemma}
  \label{lemma_lower_bound_k}
  For the number of edges $k$ of a connected $4$-regular matchstick graph in the plane we have $k\ge 11$.
\end{lemma}
  \textsc{Proof.}
  Let us assume $k=10$. Due to Equation~(\ref{eq_A_i_sum}), Lemma~\ref{lemma_inner_configuration_1},
  Lemma~\ref{lemma_area_odd}, and Table~\ref{table_max_area} we have the following possible nonzero values for the $f_i$:
  \begin{enumerate}
    \item[(a)] $f_3=14$, $f_4=t\in\mathbb{N}$, $f_{10}=1$,
    \item[(b)] $f_3=15$, $f_4=t\in\mathbb{N}$, $f_5=1$, $f_{10}=1$, or
    \item[(c)] $f_3=16$, $f_4=t\in\mathbb{N}$, $f_6=1$, $f_{10}=1$.
  \end{enumerate}
  Due to an angle sum of $(10-2)\pi$ at most $3\cdot 10-7=23$ of the inner angles of the outer face can be part
  of triangles. In a similar manner we conclude that at most $3s-3$ outer angles of an inner $s$-gon can be part
  of triangles. Since $15\cdot 3-23 -(3\cdot 5-3)=10>0$ and $16\cdot 3-23 -(3\cdot 6-3)=10>0$ there must exist one of the
  configurations of Figure~\ref{fig_configuration_1} in the cases (b) and (c), which is a contradiction to the area argument.
  
  The number of outer and inner angles of one of the configurations in Figure~\ref{fig_configuration_1}, which can be part
  of triangles, is at most $15$. Thus we can conclude from $14\cdot 3-23-1\cdot 15=4>0$  that in case (a) there must exist
  at least two subgraphs as in Figure~\ref{fig_configuration_1}, which contradicts the area argument.\hfill{\bbox}

\begin{theorem}
  Every $4$-regular matchstick graph in the plane contains at least $20$~vertices.
\end{theorem}
  \textsc{Proof.}
  Without loss of generality we can assume that the graph is connected. If the number of edges $k$ of the unbounded face
  is at least $12$, we can use Lemma~\ref{lemma_bound_4_regular} and $l\le k-3$ to conclude $n\ge 20$.
  
  If $k=11$ and $l\le k-5$ then Lemma~\ref{lemma_bound_4_regular} gives $n\ge 20$. In the remaining cases we have $k=11$ and
  $l\in\{k-4,k-3\}$ so that the area of the unbounded face is at most $\max\Big(\left(\frac{11}{4}\right)^2,
  \left(\frac{11}{3}\right)^2\cdot\frac{\sqrt{3}}{4}\Big)=\frac{121}{16}<18\cdot\frac{\sqrt{3}}{4}$. Due to
  Equation~(\ref{eq_A_i_sum}), Lemma~\ref{lemma_inner_configuration_1}, Lemma~\ref{lemma_area_odd}, and
  Table~\ref{table_max_area} we have the following possible nonzero values for the $f_i$ in this case:
  \begin{enumerate}
    \item[(a)] $f_3=15$, $f_4=t\in\mathbb{N}$, $f_{11}=1$,
    \item[(b)] $f_3=16$, $f_4=t\in\mathbb{N}$, $f_5=1$, $f_{11}=1$, or
    \item[(c)] $f_3=17$, $f_4=t\in\mathbb{N}$, $f_6=1$, $f_{11}=1$.
  \end{enumerate}
  Using similar arguments as in the proof of Lemma~\ref{lemma_lower_bound_k} we obtain a contradiction in all
  three cases.\hfill{\bbox}


\end{document}